\documentclass[12pt]{article}
\usepackage{mathrsfs}
\usepackage{t1enc}
\usepackage[latin1]{inputenc}
\usepackage[english]{babel}
\usepackage{amsmath,amsthm,amssymb}
\usepackage{amsfonts}
\usepackage{latexsym}
\usepackage{enumerate}

\usepackage{graphicx}
\usepackage[natural]{xcolor}
\usepackage{tikz}
\usepackage{tikz-3dplot}
\usetikzlibrary{shapes}
\usetikzlibrary{arrows,decorations.pathmorphing,backgrounds,positioning,fit,petri,automata}
\usetikzlibrary{positioning}

\theoremstyle{plain}
\newtheorem{thm}{Theorem}
\newtheorem{lem}[thm]{Lemma}

\newtheorem{conj}[thm]{Conjecture}

\usepackage{cases}
\theoremstyle{plain}

\theoremstyle{plain}

\theoremstyle{plain}

\usepackage[
pdfauthor={FNS},
pdftitle={Antimagic Orientation of Subdivided Caterpillars},
pdfstartview=XYZ,
bookmarks=true,
colorlinks=true,
linkcolor=blue,
urlcolor=blue,
citecolor=blue,
bookmarks=true,
linktocpage=true,
hyperindex=true
]{hyperref}

\newcommand{\arxiv}[1]{\href{http://arxiv.org/abs/#1}{\texttt{arXiv:#1}}}

\title{Antimagic orientation of subdivided caterpillars}

\author{Jessica Ferraro, Genevieve Newkirk, Songling Shan \\
{\small Department of Mathematics, Illinois State University, Normal, IL 61790, USA}}
\date{\today}

\begin{document}

\maketitle

\begin{abstract}

Let $m\ge 1$ be an integer and $G$ be a graph with $m$ edges. We say that $G$ has an antimagic orientation if $G$ has an orientation $D$ and a bijection $\tau:A(D)\rightarrow \{1,2,\ldots,m\}$ such that no two vertices in $D$ have the same vertex-sum under $\tau$, where the vertex-sum of a vertex $v$ in $D$ under $\tau$ is the sum of labels of all arcs entering $v$ minus the sum of labels of all arcs leaving $v$. Hefetz, M\"{u}tze and Schwartz [J. Graph Theory, 64: 219-232, 2010] conjectured that every connected graph admits an antimagic orientation. The conjecture was confirmed for certain classes of graphs such as regular graphs, graphs with minimum degree  at least 33, bipartite graphs with no vertex of degree zero or two, and trees including caterpillars and complete $k$-ary trees. We prove that every subdivided caterpillar admits an antimagic orientation, where a subdivided caterpillar is a subdivision of a caterpillar  $T$ such that the edges of $T$ that are not on the central path of $T$ are subdivided the same number of times. 

\smallskip 

\end{abstract}

\section{Introduction}

For two integers $p$ and $q$, let $[p,q]=\{i\in \mathbb{Z}:p\leq i
\leq q\}$. 
Let $G$ be a graph with $m\ge 1$ edges.   An \emph{antimagic labeling} of $G$ is a bijection $\tau:E(G) \rightarrow  [1,m]$ such that for any two distinct vertices $u$ and $v$, the \emph{vertex-sum}, the sum of the labels of the edges incident with a given vertex, of $u$ is distinct from the vertex-sum of $v$. A graph is said to be \emph{antimagic} if it admits an antimagic labeling.

The idea of antimagic labeling was introduced in 1990 by Hartsfield and Ringel~\cite{HR1990} and they conjectured that every connected graph and every tree other than $K_2$ is antimagic. One of the groups of researchers to continue this investigation is Kaplan, Lev and Roditty~\cite{KLR2009}, who proved that any tree having more than two vertices and at most one vertex of degree two is antimagic (also see~\cite{LWZ2014}). Other findings include proof of antimagic labelings for all regular graphs \cite{BBV2015, CLPZ2016}, graphs with average degree greater than or equal to $d_0$ for some constant $d_0$ with no isolated edge and at most one isolated vertex~\cite{E2016}, and all complete multipartite graphs other than $K_2$~\cite{AKLRY2004}. 

In 2010 Hefetz, M\"{u}tze and Schwartz introduced a variation of antimagic labeling, specifically the labeling of  digraphs~\cite{HMS2010}. Let $D$ be a digraph. We denote the vertex set and the arc set of $D$ by $V(D)$ and $A(D)$, respectively. Let $|A(D)|=m$. An antimagic labeling of $D$ is a bijection $\tau:A(D)\rightarrow  [1,m]$ such that no two vertices receive the same oriented vertex-sum, where the \emph{oriented vertex-sum} of a vertex $u\in V(D)$ is the sum of labels of all arcs entering $u$ minus the sum of labels of all arcs leaving $u$. We use $s_{[D,\tau]}(u)$ to denote the oriented vertex-sum of the vertex $u\in V(D)$ under the bijection $\tau$. For simplicity we refer to the oriented vertex-sum as the vertex-sum in the remainder of this paper as we are exclusively working with antimagic orientations. We say a graph $G$ admits an \emph{antimagic orientation} if it has an orientation $D$ such that $D$ has an antimagic labeling. Hefetz, M\"{u}tze and Schwartz~\cite{HMS2010} proposed the following conjecture.

\begin{conj} \label{orientation-conj} 
	Every connected graph admits an antimagic orientation.
\end{conj}

Note that all antimagic
bipartite graphs admit an antimagic orientation where all edges are oriented
in the same direction between the partite sets.  
In~\cite{HMS2010}, Hefetz, M\"{u}tze and Schwartz proved  Conjecture~\ref{orientation-conj}  for some classes of graphs, such as  stars, wheels, and graphs of order $n$ with minimum degree at least $c\log n$ for an absolute constant $c$. In the process the authors proved a stronger result that every orientation of these graphs is antimagic as well. Additional cases for this conjecture that have been proved already include regular graphs~\cite{HMS2010,LSWYZ2019,Y2019,SH2019}, biregular bipartite graphs with minimum degree at least two~\cite{SY2017}, Halin graphs~\cite{YCZ2019}, graphs with large maximum degree~\cite{YCOP2019}, graphs with minimum degree  at least 33 and bipartite graphs with no vertex of degree 0 or 2~\cite{S2020}. Researchers have taken particular interest in investigating trees, as we do. 
For antimagic orientation,  it is proved that Conjectur~\ref{orientation-conj} is true for caterpillars~\cite{LMS2021} (it was actually proved that caterpillars are antimagic),  complete $k$-ary trees~\cite{SH2019tree}, and lobsters~\cite{GS2020}. A \emph{caterpillar} is a tree of order at least 3 such that the removal of it's leaves produces a path. We will call a longest path in a caterpillar  the \emph{spine} and denote it by $P$ (it is easy to see that all the vertices of the caterpillar not contained in $P$ are leaves of the caterpillar), and call all the edges incident with an internal vertex of $P$ a \emph{leg} of the caterpillar. A \emph{subdivided caterpillar} $T^*$ is a subdivision of a caterpillar $T$ such that all the legs of the caterpillar are subdivided the same number of times. We again call the corresponding subdivision of 
the spine of $T$ the spine of $T^*$, and call the corresponding subdivisions of the legs of $T$ the legs of $T^*$.  

It was proved in \cite{S2020} that every bipartite graph without vertices of degree 0 or 2 admits an antimagic orientation. It suggests that constructing antimagic orientations for graphs with many vertices of degree 2 is very difficult in general. In this paper, we confirm Conjecture \ref{orientation-conj} for subdivided caterpillars whose most vertices are of degree 2.\\

\begin{thm}\label{Thm} 
	Every subdivided caterpillar admits an antimagic orientation.
\end{thm}

The remainder of the paper is organized as follows: in next section, 
we prepare some notation and preliminaries, and in Section 3, we prove Theorem~\ref{Thm}. 

\section{Notation and preliminary lemmas}

Let $T$ be a subdivided caterpillar with $m$ edges for some integer $m\ge 2$. 
Throughout the remainder of this paper, we will denote by $P$
the spine of $T$. Let 
$$
p=|E(P)|, \quad s=\text{the total number of legs of $T$}, \quad \text{and} \quad \text{$k=$ the length of each leg}. 
$$
Thus we have 
\begin{equation}\label{eq1}
m=p+ks.
\end{equation}
As it was already proved that every caterpillar admits an antimagic orientation~\cite{LMS2021}, we will 
assume from now on that 
\begin{equation}\label{eq2}
s\ge 1 \quad \text{and}\quad k\ge 2. 
\end{equation}

Furthermore, we let $P=v_0v_1\ldots v_p$ and $U=\{v_{h_1},v_{h_2},\cdots,v_{h_t}\}\subseteq V(P)$ be the set of vertices of degree at least 3 in $T$, where $h_1<h_2<\cdots<h_t$.  
We call each vertex in $U$ a \emph{joint} of $T$, and 
a joint with degree 3 in $T$ a \emph{small joint} and a \emph{big joint} otherwise. 
A leg with an endvertex as a joint is \emph{attaching} at
the joint. A leg attaching at a big joint 
is a \emph{big leg}, and a leg attaching at a small joint 
is a \emph{small leg}.  Note that $t\le s$ and when $t=s$,
then all the $t$ joints are small joints. 
%and a \emph{big joint} is a joint with degree at least 4 in $T$. 
%Likewise, a \emph{small joint} is a joint with degree 3 in $T$. 
Let $L_1, \ldots, L_s$ be the $s$ legs of $T$ and we let 
$$L_i=x_{i0}x_{i1}\ldots x_{ik},$$   for each $i\in [1,s]$, where $x_{i0} \in U$.  
Furthermore, we may assume that these legs are ordered  in consistent
with  the  joints of $T$ along $P$ from $v_0$ to $v_p$: for two joints $v_{h_i}$ and $v_{h_j}$ with $i<j$, the indices of legs attaching at $v_{h_i}$ are smaller than the indices of legs attaching at $v_{h_j}$. 

\begin{figure}[h]
	\centering

	\begin{tikzpicture}[scale=1.25, every circle node/.style={fill, inner sep=0.65ex}]
	
	%Places for the vertices
	\coordinate (P1) at (-2,0);
	\coordinate (P2) at (-1,0);
	\coordinate (P3) at (0,0);
	\coordinate (P4) at (1,0);
	\coordinate (P5) at (2,0);
	\coordinate (P6) at (3,0);
	\coordinate (P7) at (4,0);
	\coordinate (L11) at (-.5,-1);
	\coordinate (L12) at (-.5,-2);
	\coordinate (L13) at (-.5,-3);
	\coordinate (L21) at (.5, -1);
	\coordinate (L22) at (.5, -2);
	\coordinate (L23) at (.5, -3);
	\coordinate (L31) at (3,-1);
	\coordinate (L32) at (3,-2);
	\coordinate (L33) at (3,-3);
	\coordinate (EE) at (0,1);

	%For the edge labels
	\coordinate (E1) at (-1.5,0);
	\coordinate (E2) at (-.5,0);
	\coordinate (E3) at (.5,0);
	\coordinate (E4) at (1.5,0);
	\coordinate (E5) at (2.5,0);
	\coordinate (E6) at (3.5,0);
	\coordinate (E7) at (-.25,-.5);
	\coordinate (E8) at (-.5, -1.5);
	\coordinate (E9) at (-.5,-2.5);
	\coordinate (E10) at (.25,-.5);
	\coordinate (E11) at (.5, -1.5);
	\coordinate (E12) at (.5,-2.5);
	\coordinate (E13) at (3, -.5);
	\coordinate (E14) at (3, -1.5);
	\coordinate (E15) at (3, -2.5);
	
	% Drawing the lines
	\draw[thick, black] (P1) to (P2);
	\draw[thick, black,->] (P1) to (E1);
	\draw[thick, black] (P2) to (P3);
	\draw[thick, black,->] (P2) to (E2);
	\draw[thick, black] (P3) to (P4);
	\draw[thick, black,->] (P3) to (E3);
	\draw[thick, black] (P4) to (P5);
	\draw[thick, black,->] (P4) to (E4);
	\draw[thick, black] (P5) to (P6);
	\draw[thick, black,->] (P5) to (E5);
	\draw[thick, black] (P6) to (P7);
	\draw[thick, black,->] (P7) to (E6);
	\draw[thick, black] (P3) to (L11);
	\draw[thick, black,->] (L11) to (E7);
	\draw[thick, black] (L11) to (L12);
	\draw[thick, black,->] (L11) to (E8);
	\draw[thick, black] (L12) to (L13);
	\draw[thick, black,->] (L13) to (E9);
	\draw[thick, black] (P3) to (L21);
	\draw[thick, black,->] (L21) to (E10);
	\draw[thick, black] (L21) to (L22);
	\draw[thick, black,->] (L21) to (E11);
	\draw[thick, black] (L22) to (L23);
	\draw[thick, black,->] (L23) to (E12);
	\draw[thick, black] (P6) to (L31);
	\draw[thick, black,->] (L31) to (E13);
	\draw[thick, black] (L31) to (L32);
	\draw[thick, black,->] (L31) to (E14);
	\draw[thick, black] (L32) to (L33);
	\draw[thick, black,->] (L33) to (E15);

	%drawing vertices
	\node [circle, label={above:$v_0$}] at (P1) {};
	\node [circle, label={above:$v_1$}] at (P2) {};
	\node [circle, label={above,align=right:$v_2=x_{10}$\\$=x_{20}$}] at (P3) {};
	\node [circle, label={above:$v_3$}] at (P4) {};
	\node [circle, label={above:$v_4$}] at (P5) {};
	\node [circle, label={above:$v_5=x_{30}$}] at (P6) {};
	\node [circle, label={above:$v_6$}] at (P7) {};
	\node [circle, label={left:$x_{11}$}] at (L11) {};
	\node [circle, label={left:$x_{12}$}] at (L12) {};
	\node [circle, label={left:$x_{13}$}] at (L13) {};
	\node [circle, label={right:$x_{21}$}] at (L21) {};
	\node [circle, label={right:$x_{22}$}] at (L22) {};
	\node [circle, label={right:$x_{23}$}] at (L23) {};
	\node [circle, label={right:$x_{31}$}] at (L31) {};
	\node [circle, label={right:$x_{32}$}] at (L32) {};
	\node [circle, label={right:$x_{33}$}] at (L33) {};

	\end{tikzpicture} 
	\caption{Notation for $T$ when $p=6$, $s=3$, and $k=3$.}
	\label{f1}
	
\end{figure}
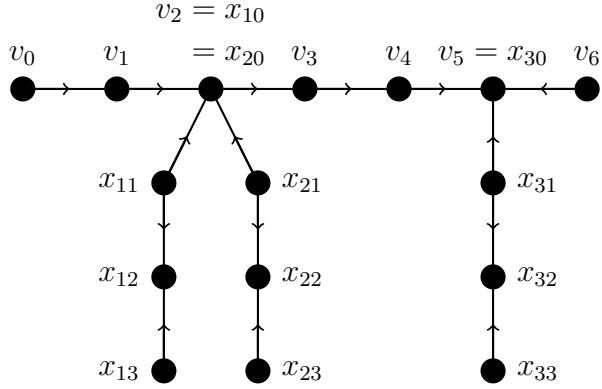

Let $L_i$ be any leg of $T$ for some $i\in [1,s]$.  In all the proofs later, $L_i$ is oriented using the following 
pattern:  
\begin{equation}\label{leg-orient}
\text{For each odd $j$ with $j\in [1,k]$, $x_{i(j-1)} \leftarrow x_{ij} $ and 
	$x_{ij} \rightarrow x_{i(j+1)}$,
}
\end{equation}
where the second arrow is defined only if  $j\le k-1$. 
See Figure~\ref{f1} below 
for an illustration of  $T$, notation defined on $T$, and the orientation of the legs of $T$.   

We will need the following result to orient and label the spine  $P$ of $T$.

\begin{lem}[\cite{GS2020}]\label{ShanLemma}  Let $P=v_0v_1\cdots v_p$ be a path and $U=\{v_{h_1}, v_{h_2},\cdots,v_{h_t}\}\subseteq V(P)$ with $|U|\ge 1$, where $p\ge 2$ and $h_1<h_2<\cdots< h_t$. Then $P$ has an orientation $\overrightarrow{P}$ and a labeling $\tau:A(\overrightarrow{P}) \rightarrow [1,p]$ such that
\begin{enumerate}[$(i)$]
	\item $s_{[\overrightarrow{P},\tau]}(v)\geq1$ for any vertex $v\in U$, and 
	    $s_{[\overrightarrow{P},\tau]}(v) \le p-1$ for every $v\in U$ with $v\ne v_{h_t}$ and $s_{[\overrightarrow{P},\tau]}(v_{h_t})  \ge 3$; 
	\item  $1\leq |s_{[\overrightarrow{P},\tau]}(v)|\leq p$ for every vertex $v\in V(P)\setminus U$; and
	\item $s_{[\overrightarrow{P},\tau]}(u)\neq s_{[\overrightarrow{P},\tau]}(v)$ for any two distinct vertices $u,v\in V(P)\setminus U$.
\end{enumerate}
\end{lem}

The second part of Lemma~\ref{ShanLemma} (i) was not specified in~\cite{GS2020}.  However, 
the conclusion is  a direct consequence of the orientation of $P$: for every $v\in U\setminus \{v_{h_t}\}$ and the two edges incident with $v$, one of the edge 
is entering $v$
and the other is  leaving $v$ in $\overrightarrow{P}$; and for the vertex $v_{h_t}$, both the edges incident with $v_{h_t}$ are entering $v_{h_t}$
in  $\overrightarrow{P}$. 

We first construct an antimahgic orientation for $T$ when $s=1$. 
\begin{lem}\label{1Leg} 
Each subdivided caterpillar with a single leg admits an antimagic orientation. 
\end{lem}

\begin{proof}

Let $T$ be a subdivided caterpillar. 
We adopt the notation introduced earlier. Thus $s=1$, $k\ge 2$ by~\eqref{eq2}, and $m=p+k$ by~\eqref{eq1}. 
We define an antimagic orientation of $T$ in two steps.
\begin{enumerate}[Step 1]
	\item Orient the edges in $E(P)$ and label them.
	
	By Lemma~\ref{ShanLemma}, $P$ has an orientation $\overrightarrow{P}$ and a labeling $\tau_1$ using numbers $1,2,\ldots,p$ satisfying the three properties described in Lemma~\ref{ShanLemma}. 
	
	\item Orient and label the leg of $T$. 
	
	Recall $L_1=x_{10} \ldots x_{1k}$ is the leg of $T$ and  $x_{10}=v_{h_1}$. 
	We orient $L_1$ using the pattern described in~\eqref{leg-orient} and 
	define $\tau_2: E(L_1) \rightarrow [p+1,m]$ according to two different cases below:  
	 $$
	\begin{cases}
	\tau_2(x_{1j}x_{1(j+1)})=m-j, \quad j\in [0,k-1], & \text{if $s_{[\overrightarrow{P},\tau_1]}(v_{h_1}) \ge m-2$},\\
	\tau_2(x_{10}x_{11})=m_1, \quad \tau_2(x_{11}x_{12})=m_2,  \quad \text{and} &\\
	 \qquad \tau_2(x_{1j}x_{1(j+1)})=m-j, \quad j\in [2,k-1], & \text{if $s_{[\overrightarrow{P},\tau_1]}(v_{h_1}) < m-2$},
	 
	\end{cases}
	$$
	where $m_1\in \{m,m-1\}$ is the number such that $m_1+s_{[\overrightarrow{P},\tau_1]}(v_{h_1})$ is even, and $\{m_2\}=\{m,m-1\}\setminus \{m_1\}$. 
\end{enumerate}

Let $T^*$ be the orientation of $T$ obtained through the two steps above, and let $$\tau: A(T^*) \rightarrow [1,m]
\,\text{such that}\,  
\tau(e)=
\begin{cases}
\tau_1(e)  & \text{if $e\in A(\overrightarrow{P})$},\\
\tau_2(e) & \text{if $e\in A(T^*)\setminus A(\overrightarrow{P})$}, 
\end{cases}
$$
where note that we treat an arc $e$ of $T^*$ and the corresponding 
un-directed $e$ as the same element in defining $\tau$. (This convention 
will be used throughout the paper when we defining an antimagic labeling 
of an orientation of a graph.) 

 We show next that $\tau$ is an antimagic labeling of $T^*$.  Let $u,v\in V(T^*)$ be any two distinct 
vertices. By the labeling, it is clear that $s_{[T^*,\tau]}(u) \ne s_{[T^*,\tau]}(v)$
if $u,v\in V(P)$ or $u,v\in V(L_1)\setminus \{v_{h_1}\}$. 
If $u\in V(P)\setminus\{v_{h_1}\}$ and $v\in V(L_1)\setminus \{v_{h_1}\}$, then $|s_{[T^*,\tau]}(u)|\in [1,p]$
and $|s_{[T^*,\tau]}(v)|\in [p+1,2m-1]$, and so they are different. 
Thus we assume that $u=v_{h_1}$ and $v\in V(L_1)\setminus \{v_{h_1}\}$.
If $s_{[\overrightarrow{P},\tau_1]}(v_{h_1}) \ge m-2$, then $s_{[T^*,\tau]}(u) \ge m-2+\tau_2(x_{10}x_{11})=m-2+m= 2m-2$. 
Note that $s_{[T^*,\tau]}(x_{1j})=(-1)^j(2m-2j+1)$ for each $j\in [1,k-1]$
and $|s_{[T^*,\tau]}(x_{1k})|=p+1$. Thus 
for even $j\in [1,k]$, $s_{[T^*,\tau]}(x_{1j}) \le 2m-3$. Hence $s_{[T^*,\tau]}(v_{h_1}) \ne s_{[T^*,\tau]}(v)$. 
Therefore we assume  $s_{[\overrightarrow{P},\tau_1]}(v_{h_1}) < m-2$. 
Then by  Lemma~\ref{ShanLemma} (i) and the definition of $\tau_2$, we know that $m+2=3+m-1 \le s_{[T^*,\tau]}(v_{h_1}) \le m-3+m=2m-3$ and  $s_{[T^*,\tau]}(v_{h_1})$ 
is  even. 
Since $s_{[T^*,\tau]}(x_{12}) \in\{2m-2,2m-3\}$,  $s_{[T^*,\tau]}(x_{1j})=(-1)^j(2m-2j+1)$ is odd for $j\in [3,k-1]$, 
and $s_{[T^*,\tau]}(x_{1k}) \le m$, 
we again have $s_{[T^*,\tau]}(v_{h_1}) \ne s_{[T^*,\tau]}(v)$. 
This finishes the proof. 
%Given the structure of $T$, $m = p+k$. Given the use of this pattern and the structure of $T$, the absolute value of the sums of the degree 1 vertices fall into the interval $[1,p+1]$. Considering the degree 2 vertices, they will either fall in the path or the leg. When in the path the absolute value of their sums will be between 1 and $p$. When they fall in the leg, there sum will be $(-1)^d[2m-2d+1]$, with one exception in Case 2 and Case 3. These typical sums will always be odd. In the exception, when $\partial(v_{h_t})$ is odd in Case 2 and $\partial(v_{h_t})$ is even in Case 3, consider the vertex $x_2$. For this vertex the sum will be $2m-2$, which will result in an even value. Though the sum of the degree 3 vertex is also even, $s(v_{h_t})<m-2+m=2m-2$. Therefore, the two sums will never equal though they have the same parity. This means that $T$ admits an antimagic orientation when $s=1$.\\
\end{proof}

\section{Proof of Theorem 2}

\begin{proof}
	Let $T$ be a subdivided caterpillar. 
	We adopt the notation introduced in Section 2. 
By Lemma \ref{1Leg} and \eqref{eq2}, we assume $k \geq 2$ and $s\geq2$.  
We construct an antimagic orientation of $T$
below. 

\smallskip

{\bf Step 1: Orient the edges in $E(P)$ and label them.}

\smallskip 

By Lemma~\ref{ShanLemma}, $P$ has an orientation $\overrightarrow{P}$ and a labeling $\tau_1$ using numbers $1,2,\ldots,p$ satisfying  the three properties described in Lemma~\ref{ShanLemma}.

\smallskip 

{\bf Step 2: Orient and label the legs of $T$.}
\smallskip 

We orient each  leg of $T$ using the pattern described in~\eqref{leg-orient}. 
Together with the orientation of $P$ in Step 1, we have obtained an 
orientation $\overrightarrow{T}$ of $T$. 

Next, we will 
assign labels in $[p+1,m]$ to the edges contained in legs of $T$. Define 
$$E_1=\{x_{i0}x_{i1}: i\in [1,s]\},$$
to be the set of the edges from the legs that are incident with 
the joints of $T$. 
\smallskip 

{\bf Step 2.1: Assign  labels in $[m-s+1,m]$ to  edges in $E_1$.}
\smallskip

We let $M \subseteq E_1$ be a matching of $T$ with size $|U|=t$ and saturating $U$, recall $U=\{v_{h_1}, \ldots, v_{h_t}\}$
is the set of vertices of $T$ of degree at least 3 in $T$.  We arbitrarily assign
labels in $[m-(s-t)+1,m]$ to edges in $E_1\setminus M$ such that distinct edge receive 
distinct label. 
Denote by $\tau_2^*$  the current labeling of $T^*_1$, where 
$T^*_1$ consists of  $\overrightarrow{P}$ and the oriented 
edges from $E_1\setminus M$.  
Now for each  vertex $v\in U$, we compute  $s_{[T^*_1,\tau_2^*]}(v)$, and assume that 
$$
s_{[T^*_1,\tau_2^*]}(x_1) \ge \ldots  \ge s_{[T^*_1,\tau_2^*]}(x_t),  
$$
where $\{x_1,\ldots, x_t\}$ is a permutation of the vertices of $U$. 
Let $T_1$ be the union of $T^*_1$ and those edges from $M$ together with their orientation. 
Now define $$\tau_2: A(T_1) \rightarrow [m-s+1,m]  \, \text{with} \, 
\tau_2(e)=
\begin{cases}
\tau_2^*(e),   & \text{if $e\in A(T_1^*)$,}\\
m-(s-t)+1-i, & \text{if $e\in M$ and is incident with $x_i$}. 
\end{cases}
 $$
 By this definition of $\tau_2$, we have 
 \begin{equation}\label{U-vertex}
s_{[T_1, \tau_2]}(x_1) >\ldots > s_{[T_1, \tau_2]}(x_t). 
 \end{equation}
 
 \smallskip 
 
 {\bf Step 2.2: Assign  labels in $[p+1,m]$ to edges in $E(T)\setminus E(P)$.}
 \smallskip 

Let $i\in [1,s]$.
For the leg $L_i$, assume $s_{[T_1,\tau_2]}(x_{i0}x_{i1})=a_i$.
 Then the labels will be used for edges of $L_i$
will be the set 
$$
A_i=\{a_i, a_i-s, a_i-2s, \ldots, a_i-(k-1)s \}. 
$$
It is clear that  $\bigcup_{i=1}^s A_i =[p+1, m]$. 
Let $f:  E(T)\setminus E(P)\rightarrow [p+1,m]$ be a bijection. For $e\in E(T)\setminus E(P)$, assume $e\in E(L_i)$ for some $i\in [1,s]$.
 We  define $f$ 
according to three different cases (case 2 and case 3 can be combined, but we separate them for clarity), see Figure~\ref{f2} for an illustration.

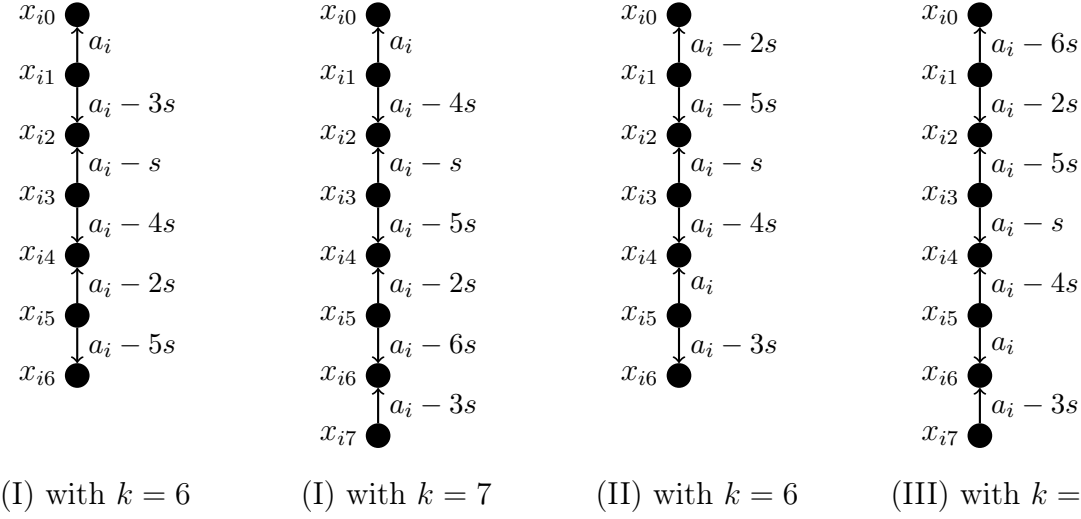
\begin{figure}[!htb]
	\begin{center}
		\begin{tikzpicture}[scale=0.8]

		{\tikzstyle{every node}=[draw ,circle,fill=black, minimum size=0.3cm,
			inner sep=0pt]
			\draw[black,thick](1,1) node[label={left: $x_{i0}$\,}] (x0)  {};
			\draw[black,thick](1,0) node[label={left: $x_{i1}$\,}] (x1)  {};
			\draw[black,thick](1,-1) node[label={left: $x_{i2}$\,}] (x2)  {};
			\draw[black,thick](1,-2) node[label={left: $x_{i3}$\,}] (x3)  {};
			\draw[black,thick](1,-3) node[label={left: $x_{i4}$\,}] (x4)  {};
			\draw[black,thick](1,-4) node[label={left: $x_{i5}$\,}] (x5)  {};
			\draw[black,thick](1,-5) node[label={left: $x_{i6}$\,}] (x6)  {};
		}
	
		\path[draw,thick,black,->]
		(x1) edge node[name=la,pos=0.5, right] {\small $a_i$} (x0) 
		(x1) edge node[name=la,pos=0.5, right] {\small $a_i-3s$} (x2)
		(x3) edge node[name=la,pos=0.5, right] {\small $a_i-s$} (x2)
		(x3) edge node[name=la,pos=0.5, right] {\small $a_i-4s$} (x4)
		(x5) edge node[name=la,pos=0.5, right] {\small $a_i-2s$} (x4)
		(x5) edge node[name=la,pos=0.5, right] {\small $a_i-5s$} (x6); 
		
		\node at (1.3,-7) {  (I) with $k=6$};
	
	\begin{scope}[shift={(5,0)}]
	{\tikzstyle{every node}=[draw ,circle,fill=black, minimum size=0.3cm,
		inner sep=0pt]
		\draw[black,thick](1,1) node[label={left: $x_{i0}$\,}] (x0)  {};
		\draw[black,thick](1,0) node[label={left: $x_{i1}$\,}] (x1)  {};
		\draw[black,thick](1,-1) node[label={left: $x_{i2}$\,}] (x2)  {};
		\draw[black,thick](1,-2) node[label={left: $x_{i3}$\,}] (x3)  {};
		\draw[black,thick](1,-3) node[label={left: $x_{i4}$\,}] (x4)  {};
		\draw[black,thick](1,-4) node[label={left: $x_{i5}$\,}] (x5)  {};
		\draw[black,thick](1,-5) node[label={left: $x_{i6}$\,}] (x6)  {};
		\draw[black,thick](1,-6) node[label={left: $x_{i7}$\,}] (x7)  {};
	}
	
	\path[draw,thick,black,->]
	(x1) edge node[name=la,pos=0.5, right] {\small $a_i$} (x0) 
	(x1) edge node[name=la,pos=0.5, right] {\small $a_i-4s$} (x2)
	(x3) edge node[name=la,pos=0.5, right] {\small $a_i-s$} (x2)
	(x3) edge node[name=la,pos=0.5, right] {\small $a_i-5s$} (x4)
	(x5) edge node[name=la,pos=0.5, right] {\small $a_i-2s$} (x4)
	(x5) edge node[name=la,pos=0.5, right] {\small $a_i-6s$} (x6)
	(x7) edge node[name=la,pos=0.5, right] {\small $a_i-3s$} (x6); 
	
	\node at (1.3,-7) {  (I) with $k=7$};
	
	\end{scope}
		
		\begin{scope}[shift={(10,0)}]
			{\tikzstyle{every node}=[draw ,circle,fill=black, minimum size=0.3cm,
			inner sep=0pt]
			\draw[black,thick](1,1) node[label={left: $x_{i0}$\,}] (x0)  {};
			\draw[black,thick](1,0) node[label={left: $x_{i1}$\,}] (x1)  {};
			\draw[black,thick](1,-1) node[label={left: $x_{i2}$\,}] (x2)  {};
			\draw[black,thick](1,-2) node[label={left: $x_{i3}$\,}] (x3)  {};
			\draw[black,thick](1,-3) node[label={left: $x_{i4}$\,}] (x4)  {};
			\draw[black,thick](1,-4) node[label={left: $x_{i5}$\,}] (x5)  {};
			\draw[black,thick](1,-5) node[label={left: $x_{i6}$\,}] (x6)  {};
		}
		
		\path[draw,thick,black,->]
		(x1) edge node[name=la,pos=0.5, right] {\small $a_i-2s$} (x0) 
		(x1) edge node[name=la,pos=0.5, right] {\small $a_i-5s$} (x2)
		(x3) edge node[name=la,pos=0.5, right] {\small $a_i-s$} (x2)
		(x3) edge node[name=la,pos=0.5, right] {\small $a_i-4s$} (x4)
		(x5) edge node[name=la,pos=0.5, right] {\small $a_i$} (x4)
		(x5) edge node[name=la,pos=0.5, right] {\small $a_i-3s$} (x6); 
		
		\node at (1.3,-7) {  (II) with $k=6$};
		
		\end{scope}

			\begin{scope}[shift={(15,0)}]
		{\tikzstyle{every node}=[draw ,circle,fill=black, minimum size=0.3cm,
			inner sep=0pt]
			\draw[black,thick](1,1) node[label={left: $x_{i0}$\,}] (x0)  {};
			\draw[black,thick](1,0) node[label={left: $x_{i1}$\,}] (x1)  {};
			\draw[black,thick](1,-1) node[label={left: $x_{i2}$\,}] (x2)  {};
			\draw[black,thick](1,-2) node[label={left: $x_{i3}$\,}] (x3)  {};
			\draw[black,thick](1,-3) node[label={left: $x_{i4}$\,}] (x4)  {};
			\draw[black,thick](1,-4) node[label={left: $x_{i5}$\,}] (x5)  {};
			\draw[black,thick](1,-5) node[label={left: $x_{i6}$\,}] (x6)  {};
			\draw[black,thick](1,-6) node[label={left: $x_{i7}$\,}] (x7)  {};
		}
		
		\path[draw,thick,black,->]
		(x1) edge node[name=la,pos=0.5, right] {\small $a_i-6s$} (x0) 
		(x1) edge node[name=la,pos=0.5, right] {\small $a_i-2s$} (x2)
		(x3) edge node[name=la,pos=0.5, right] {\small $a_i-5s$} (x2)
		(x3) edge node[name=la,pos=0.5, right] {\small $a_i-s$} (x4)
		(x5) edge node[name=la,pos=0.5, right] {\small $a_i-4s$} (x4)
		(x5) edge node[name=la,pos=0.5, right] {\small $a_i$} (x6)
		(x7) edge node[name=la,pos=0.5, right] {\small $a_i-3s$} (x6); 
		
		\node at (1.3,-7) {  (III) with $k=7$};
		
		\end{scope}

		\end{tikzpicture}
		-	  	\end{center}
	\caption{Three different labeling patterns of leg $L_i$.}
	\label{f2}
\end{figure}

\begin{enumerate}[(I)]
	\item  For even $j \in [0,k-2]$, let  
	 $$f(x_{ij}x_{i(j+1)})=a_i-\frac{j}{2}s, \quad  f(x_{i{j+1}}x_{i(j+2)})=a_i-\left(\left \lceil \frac{k}{2} \right \rceil+\frac{j}{2} \right)s. $$
	 When $k$ is odd, let  $f(x_{i(k-1)}x_{ik})=a_i-\frac{k-1}{2}s$. 
	 \item $k$ is even, for even $j \in [0,k-2]$, let  
	$$f(x_{ij}x_{i(j+1)})=a_i-(k-1)s+\left(\frac{k}{2}+\frac{j}{2} \right)s, \quad  f(x_{i{j+1}}x_{i(j+2)})=a_i-(k-1)s+\frac{j}{2}s.$$ 
	
		\item  $k$ is odd, let  $f(x_{i(k-1)}x_{ik})=a_i-\frac{k-1}{2}s$ and for  even $j \in [0,k-3]$, let 
	$$f(x_{ij}x_{i(j+1)})=a_i-(k-1)s+\frac{j}{2}s, \quad  f(x_{i{j+1}}x_{i(j+2)})=a_i-(k-1)s+\left(\frac{k+1}{2}+\frac{j}{2} \right)s. $$ 
\end{enumerate}

We define a bijection $\tau_3:  E(T)\setminus E(P)\rightarrow [p+1,m]$ as follows: for $e\in E(T)\setminus E(P)$, assume $e\in E(L_i)$ for some $i\in [1,s]$. Then 
$$
\tau_3 (e)=
\begin{cases}
f(e) \quad \text{in (I)}, & \text{if $L_i$ is a big leg}, \\
 f(e) \quad \text{in (II)}, & \text{if $k$ is even and $L_i$ is a small leg such that $x_{i0} \ne v_{h_t}$}, \\ 
   f(e) \quad \text{in (III)}, & \text{if $k$ is odd and $L_i$ is a small leg such that $x_{i0} \ne v_{h_t}$}. 
\end{cases}
$$
Assume now that $L_i$ is a small leg and $x_{i0} = v_{h_t}$.  Let 
$$
\tau_3 (e)=
\begin{cases}
f(e) \quad \text{in (I)}, & \text{if  $s_{[\overrightarrow{P}, \tau_1]}(v_{h_t}) \ge m-s$}, \\
f(e) \quad \text{in (II)}, & \text{if $k$ is even and $s_{[\overrightarrow{P}, \tau_1]}(v_{h_t})  \le p$}, \\
f(e) \quad \text{in (III)}, & \text{if $k$ is odd and $s_{[\overrightarrow{P}, \tau_1]}(v_{h_t})  \le p$}. 
\end{cases}
$$

Lastly, assume $L_i$ is a small leg, $x_{i0} = v_{h_t}$, and $ p+1\le s_{[\overrightarrow{P}, \tau_1]}(v_{h_t}) \le  m-s-1$. 
Assume, without loss of generality that $x_{i0}=v_{h_t}=x_\ell$ for some $\ell\in [1,s]$,  where recall $x_\ell$ is defined in~\eqref{U-vertex}. 
Since $s_{[\overrightarrow{P}, \tau_1]}(v_{h_t})\ge p+1$ 
and $s_{[\overrightarrow{P}, \tau_1]}(v_{h_i}) \le p$  for $v_{h_i} \in U\setminus\{v_{h_t}\}$, 
and  by the definition of $\tau_2$, it follows that if $\ell \ge 2$, then $T$
has  big legs, and if $T$ has no big  leg, then $\ell=1$. 
We define 
$$
\tau_3 (e)=
\begin{cases}
f(e) \quad \text{in (I)}, & \text{if  $ \ell \ge 2$}, \\
f(e) \quad \text{in (II)}, & \text{if $k$ is even and $ \ell=1$}, \\
f(e) \quad \text{in (III)}, & \text{if $k$ is odd and $\ell=1$}. 
\end{cases}
$$

\smallskip 

{\bf Step 2.3: Modify $\tau_3$ defined in Step 2.2 to avoid same vertex-sums.}
\smallskip 

For the bijection $\tau_3$ defined in Step 2.2, it might happen that 
one vertex-sum of some degree 2 vertices from the legs of $T$ is the same 
as the vertex-sum of the vertex $v_{h_t}$. For this reason, we will modify 
$\tau_3$ slightly so that under the modification, the vertex-sum of 
$v_{h_t}$ is not the same  as those of degree 2 vertices $x_{ij}$
for all $i\in [1,s]$ and all  $j\in [1,k-1]$. 

Recall that $v_{h_t}=x_\ell$. When $v_{h_t}$ is a small joint,  we assume $L_q$, for some $q\in [1,s]$,  
is the leg of $T$ with  $x_{q0}=v_{h_t}$.  If $\ell \ge 2$ and so $q\ge 2$, 
assume, by relabeling the legs if necessary, that $L_{q-1}$
is the leg of $T$ with  $x_{(q-1)0}=x_{\ell-1}$ and $x_{(q-1)0}x_{(q-1)1}\in M$. 
Thus by the definition of $\tau_2$ in Step 2.1, we have 
$$
\tau_2(x_{(q-1)0}x_{(q-1)1})=m-(s-t)+1-(\ell-1) \quad \text{and} \quad \tau_2(x_{q0}x_{q1})=m-(s-t)+1-\ell. 
$$
If $\ell =1$, then $s=t$ and $q=1$ by the way of naming the legs of $T$. 
Thus $L_{q+1}=L_2$
is the leg of $T$ with  $x_{(q+1)0}=x_{\ell+1}$ and $x_{(q+1)0}x_{(q+1)1}\in M$. 
Thus by the definition of $\tau_2$ in Step 2.1, we have 
$$
\tau_2(x_{q0}x_{q1})=m \quad \text{and} \quad \tau_2(x_{(q+1)0}x_{(q+1)1})=m-1. 
$$

We define $\tau_4: E(T)\setminus E(P)\rightarrow [p+1,m]$ by modifying $\tau_3$ 
as below.  If $v_{h_t}$ is a big joint, we let $\tau_4=\tau_3$.
Thus we assume $v_{h_t}$ is a small joint. 
If  $  s_{[\overrightarrow{P}, \tau_1]}(v_{h_t}) \le p$ or $  s_{[\overrightarrow{P}, \tau_1]}(v_{h_t})  \ge m-s$, we let 
$\tau_4=\tau_3$. Thus we assume $v_{h_t}$ is a small joint and $ p+1\le s_{[\overrightarrow{P}, \tau_1]}(v_{h_t}) \le  m-s-1$.  Under this assumption, 
we modify $\tau_3$ in two different subcases. 
\begin{enumerate}[(a)]
	\item If (a): $s$ is even and $s_{[\overrightarrow{P}, \tau_1]}(v_{h_t})+m-(s-t)+1-\ell$ is odd, (b): $s$ is odd,  $k$ is even and $s_{[\overrightarrow{P}, \tau_1]}(v_{h_t})+m-(s-t)+1-\ell$ has a different parity than $k/2-1$, or (c): $s$ is odd, $k$ is odd and $s_{[\overrightarrow{P}, \tau_1]}(v_{h_t})+m-(s-t)+1-\ell$ has a different parity than $(k-1)/2$,  then we let $\tau_4=\tau_3$. 
	\item In all other cases we let $\tau_4$ be obtained 
	from $\tau_3$ by switching the labels on the leg $L_{q-1}$ with that of $L_q$
	if $\ell \ge 2$, and by switching the labels on the leg $L_{q}$ with that of $L_{q+1}$
	if $\ell =1$. 
\end{enumerate}
 
From now on, {\bf we will call the leg $L_q$ big if it was labeled using the pattern $f$ in (I), 
and small otherwise}.

We now let $\tau: A(\overrightarrow{T}) \rightarrow [1,m]$ be the bijection obtained from $\tau_1$
and $\tau_4$ by letting $\tau(e)=\tau_1(e)$ if $e\in A(\overrightarrow{P})$ and 
$\tau(e)=\tau_4(e)$ if $e\in A(\overrightarrow{T})\setminus A(\overrightarrow{P})$.   For notation simplicity, we write $s_{[\overrightarrow{T},\tau]}$ as $s_{\tau}$ in the rest of the proof. 
Let $u,v\in V(\overrightarrow{T})$ be any two distinct vertices. 
We show $s_\tau(u) \ne  s_\tau(v)$ in the following 6 cases,  which implies that $\tau$
is an antimagic labeling of $\overrightarrow{T}$.  

\begin{enumerate}
	\item []
\begin{enumerate}
	\item [Case 1:] both $u$ and $v$ are leaves of $T$. 
		\item [Case 2:] $u$ is a leaf of $T$ and $v\in U$. 
		
	\item [Case 3:] $u,v\in U$. 
		\item [Case 4:] $u$ is a leaf of $T$ and $v$ is a degree 2 vertex of $T$. 
			\item [Case 5:] both $u$ and $v$ are degree 2 vertices of $T$. 
	\item [Case 6:]  $u$ is a degree 2 vertex of $T$ and $v \in U$.

\end{enumerate}
\end{enumerate}

{\bf \noindent  Case 1}: both $u$ and $v$ are leaves of $T$. \quad 

 As $\tau$ is a bijection and so distinct edges receive distinct labels under $\tau$,  it is clear that $s_\tau(u) \ne  s_\tau(v)$. 

{\bf \noindent  Case 2}: $u$ is a leaf of $T$ and $v\in U$.  \quad

Note that $s_\tau(v) \ge p+1$. Thus  $s_\tau(u) \ne  s_\tau(v)$ 
when $u\in V(P)$ 
or $k$ is odd,  since $s_\tau(u)  \le p$ or $s_\tau(u) <0$.
So assume $u\in V(T) \setminus V(P)$ and $k$ is even. 
When  $k$ is even, then  by (I) and (II) of the definition of $f$, we have  $s_\tau(u) \le m-ks/2s$. 
However,  $s_\tau(v) \ge s_{[\overrightarrow{P}, \tau_1]}(v)+(m-s+1)-(k-1)s+ks/2>m+1-ks/2>s_\tau(u)$. 
Thus again, $s_\tau(u) \ne  s_\tau(v)$.  

{\bf \noindent  Case 3}: $u,v\in U$.  \quad

Recall  $T_1$  consists of  $\overrightarrow{P}$ and the oriented 
edges from $E_1$.  
We let $\pi:A(T_1) \rightarrow [1,p]\cup \{\tau_3(x_{i0}x_{i1}): i\in [1,s]\}$ be the bijection obtained from $\tau_1$
and $\tau_3$ by letting $\pi(e)=\tau_1(e)$ if $e\in A(\overrightarrow{P})$ and 
$\pi(e)=\tau_3(e)$ if $e\in A(T_1)\setminus A(\overrightarrow{P})$.

By the definition of $\tau_3$, we have $\tau_3(e) \le \tau_2(e)$ for each $e\in E_1$. Thus by~\eqref{U-vertex},
we have 
 \begin{equation}\label{U-vertex2}
s_{[T_1, \pi]}(x_1) >\ldots > s_{[T_1, \pi]}(x_t). 
\end{equation}

Recall $\tau_4$ was obtained by modifying $\tau_3$ either in terms of $\tau_4=\tau_3$ or by exchanging the labels 
on $L_q$ with that on $L_{q-1}$ if $\ell \ge 2$
and with that on $L_{q+1}$ if $\ell=1$, where recall $x_\ell= v_{h_t}$. 
If $\tau_4=\tau_3$, then for each $i\in [1,t]$, we have  $s_{\tau}(x_i) = s_{[T_1, \pi]}(x_i)$
and thus $s_\tau(u) \ne  s_\tau(v)$ by~\eqref{U-vertex2}. 

Thus we assume $\tau_4\ne \tau_3$, which implies that $v_{h_t}$ is a small joint and $ p+1\le s_{[\overrightarrow{P}, \tau_1]}(v_{h_t}) \le  m-s-1$ by
the definition of $\tau_4$.  
When $\ell \ge 2$, we have $s_{\tau}(x_i) = s_{[T_1, \pi]}(x_i)$ for all $i\in [1,t]\setminus \{\ell-1,\ell\}$, 
$s_{\tau}(x_{\ell-1}) = s_{[T_1, \pi]}(x_{\ell-1})-1$ and $s_{\tau}(x_\ell) = s_{[T_1, \pi]}(x_\ell)+1$. 
Note that $x_{\ell-1}$ is a big joint  and $x_\ell =v_{h_t}$ is  a small joint. 
Thus $$s_{[T_1, \pi]}(x_{\ell-1}) \ge s_{[\overrightarrow{P}, \tau_1]}(x_{\ell-1})+m-s+t+1+m-(s-t)+1-(\ell-1) >2m+3-(s-t)-\ell-s$$ 
and $$s_{[T_1, \pi]}(x_{\ell}) = s_{[\overrightarrow{P}, \tau_1]}(x_{\ell})+m-(s-t)+1-\ell \le m-s-1+m-(s-t)+1-\ell=2m-(s-t)-\ell-s.$$ 
Therefore,  $s_{\tau}(x_{\ell-1})=s_{[T_1, \pi]}(x_{\ell-1})-1>s_{[T_1, \pi]}(x_\ell)+1=s_{\tau}(x_\ell)$ and  
so the strict  inequalities in~\eqref{U-vertex2} still hold with respect to $\tau$. 
As a consequence  $s_\tau(u) \ne  s_\tau(v)$. 

We then assume $\ell=1$. We have $s_{\tau}(x_i) = s_{[T_1, \pi]}(x_i)$ for all $i\in [1,t]\setminus \{1,2\}$, 
$s_{\tau}(x_{1}) = s_{[T_1, \pi]}(x_{1})-1$ and $s_{\tau}(x_2) = s_{[T_1, \pi]}(x_2)+1$. 
In this cases, both $x_1=v_{h_t}$ and $x_2$ are small joints.  Thus 
$$s_{[T_1, \pi]}(x_{1}) = s_{[\overrightarrow{P}, \tau_1]}(x_1)+m \ge p+1+m$$ 
and $$s_{[T_1, \pi]}(x_2) = s_{[\overrightarrow{P}, \tau_1]}(x_2)+m-1 \le p-1+m-1,$$ 
by Lemma~\ref{ShanLemma} (i). 
Therefore,  $s_{\tau}(x_1)=s_{[T_1, \pi]}(x_{1})-1 >s_{[T_1, \pi]}(x_2)+1=s_{\tau}(x_2)$ and so the street inequalities in~\eqref{U-vertex2} still hold with respect to $\tau$. 
This again gives  $s_\tau(u) \ne  s_\tau(v)$.

{\bf \noindent  Case 4}: $u$ is a leaf of $T$ and $v$ is a degree 2 vertex of $T$.  

 If $u, v\in V(P)$, then we have $s_\tau(u) \ne  s_\tau(v)$ by Lemma~\ref{ShanLemma}. 
 Thus we assume that $|V(P)\cap \{u,v\}| \le 1$. 
 Note that although $\tau_4$ is a modification of $\tau_3$, the 
 set of vertex-sums induced on the vertices of $V(T)\setminus V(P)$
 are the same under both of them, and the bijection $f$. 
 We here derive formulas for those vertex-sums. 
Let $i\in [1,s]$  and $L_i$ be a leg. By the definition of $f$, when $k$ is even,
we have 
\begin{equation}\label{k-even-sum} 
\begin{cases}
s_\tau(x_{ij})=(-1)^j \left(2a_i-(\frac{k}{2}+j-1)s \right), \, j\in [1,k-1], &  \text{if $L_i$ is a big leg},\\ 
s_\tau(x_{ij})= (-1)^j \left(2a_i-(2(k-1)-\frac{k}{2}-(j-1))s \right), \, j\in [1,k-1],&  \text{if $L_i$ is a  small leg}.\\
\phantom{s_\tau(x_{ij})}=(-1)^j \left(2a_i-(\frac{3k}{2}-1-j)s \right). 
\end{cases}
\end{equation} 
By the definition of $f$, when $k$ is odd,
we have 
\begin{equation}\label{k-odd-sum} 
\begin{cases}
 s_\tau(x_{ij})=(-1)^j \left(2a_i-(\frac{k-1}{2}+j)s \right), \, j\in [1,k-1],&  \text{if $L_i$ is a big leg},\\ 
s_\tau(x_{ij})=(-1)^j \left(2a_i-(2(k-1)-\frac{k+1}{2}-(j-1))s \right), \, j\in [1,k-1],&  \text{if $L_i$ is a  small leg}. \\
\phantom{s_\tau(x_{ij})}=(-1)^j \left(2a_i-(\frac{3(k-1)}{2}-j)s \right).
\end{cases}
\end{equation} 

Similarly, when $k$ is even,
we have 
\begin{equation}\label{k-even-sum2} 
\begin{cases}
s_\tau(x_{ik})=a_i-(k-1)s,  &  \text{if $L_i$ is a big leg},\\ 
s_\tau(x_{ik})= a_i- \frac{k}{2}s, &  \text{if $L_i$ is a  small leg}.\\
\end{cases}
\end{equation} 
When $k$ is odd,
we have 
\begin{equation}\label{k-odd-sum2} 
\begin{cases}
s_\tau(x_{ik})=-(a_i-\frac{k-1}{2}s),   &  \text{if $L_i$ is a big leg},\\ 
s_\tau(x_{ik})=-(a_i-\frac{k-1}{2}s),&  \text{if $L_i$ is a  small leg}. \\
\end{cases}
\end{equation} 

By the formulas above we see that 
 $|s_\tau(v)|  \ge 2(m-s+1)-3s(k-1)/2 \ge p+1$  if $d_T(v)=2$ and $v \not \in V(P)$.
 By the definition of $\tau_1$ and Lemma~\ref{ShanLemma},   we have  $|s_\tau(v)| \le p$  
 if $d_T(v)=2$ and $v  \in V(P)$. For any leaves $u$ of $T$, 
 we have $ p+1\le |s_\tau(u)| \le m-s(k-1)/2$ if $u\not\in V(P)$
 and   $|s_\tau(u)| \le p$ if $u\in V(P)$. 
 Thus we have $s_\tau(u) \ne  s_\tau(v)$ if $|V(P)\cap \{u,v\}| \le 1$. 
 Hence we assume $u,v\not\in V(P)$. 
 Then $$|s_\tau(v)|-|s_\tau(u)| \ge 2(m-s+1)-3s(k-1)/2-( m-s(k-1)/2)=m+2-s-ks>0,$$
as $m=p+ks>s+ks$.  Thus $s_\tau(u) \ne  s_\tau(v)$. 

{\bf \noindent  Case 5}: both $u$ and $v$ are degree 2 vertices of $T$.   

If $u, v\in V(P)$, then we have $s_\tau(u) \ne  s_\tau(v)$ by Lemma~\ref{ShanLemma}. 
If $u\in V(P)$ and $v\notin V(P)$, then $|s_\tau(u)| \le p <| s_\tau(v)|$. Thus 
we assume $u,v\in V(T)\setminus V(P)$ and 
$u=x_{ij}\in V(L_i)$ and $v=x_{hr} \in V(L_h)$ for some $i,h\in [1,s]$ and $j,r\in [1,k-1]$.  It is clear that $s_\tau(u) \ne  s_\tau(v)$ if $j$ and $r$ have different 
parities. Thus we assume $j \equiv r \pmod 2$. 

We have three subcases to analysis:  both $L_i$ and $L_h$ are big legs,  $L_i$ is a big leg and $L_h$ is a small leg, and  both $L_i$ and $L_h$ are small legs (recall that  the leg $L_q$ with $x_{q0}=v_{h_t}$
is called small only if it is labeled using the definition of $f$ in (II) or (III)). 
When $k$ is even, by~\eqref{k-even-sum}, 
\begin{eqnarray*}
|s_\tau(x_{ij})- s_\tau(x_{hr})|= & \begin{cases}
	|2a_i-2a_h+(r-j)s|, & \text{if both $L_i$ and $L_h$ are big legs}, \\
	|2a_i-2a_h+(k-r-j)s|, & \text{if  $L_i$ is big and $L_h$ is small}, \\
	|2a_i-2a_h+(j-r)s|, & \text{if both $L_i$ and $L_h$ are small legs}. 
\end{cases}
\end{eqnarray*} 

Since $k$ is even and $j \equiv r \pmod 2$,  we know that both $|r-j|$
and $|k-r-j|$ are even. Since $u$ and $v$
are distinct vertices, we know that if $i=h$, then $ j\ne r$
and if $j=r$ then $i\ne h$.  If $L_i$ is big and $L_h$ is small,
then $a_i\ne a_h$. 
These facts together with the fact that 
$a_i,a_h\in [m-s+1,m]$ and so $|2a_i-2a_h| \le 2(s-1)$, imply  
 $|s_\tau(x_{ij})- s_\tau(x_{hr})| \ne 0$. 
Thus $s_\tau(u) \ne  s_\tau(v)$. 

We then assume that $k$ is odd. Then  
by~\eqref{k-odd-sum}, 
\begin{eqnarray*}
	|s_\tau(x_{ij})- s_\tau(x_{hr})|= & \begin{cases}
		|2a_i-2a_h+(r-j)s|, & \text{if both $L_i$ and $L_h$ are big legs}, \\
		|2a_i-2a_h+(k-1-r-j)s|, & \text{if  $L_i$ is big and $L_h$ is small}, \\
		|2a_i-2a_h+(j-r)s|, & \text{if both $L_i$ and $L_h$ are small legs}. 
	\end{cases}
\end{eqnarray*} 
Since $k-1$ is even and $j \equiv r \pmod 2$,  we know that both $|r-j|$
and $|k-1-r-j|$ are even. Since $u$ and $v$
are distinct vertices, we know that if $i=h$, then $ j\ne r$
and if $j=r$ then $i\ne h$.  If $L_i$ is big and $L_h$ is small,
then $a_i\ne a_h$. 
These facts together with the fact that 
$a_i,a_h\in [m-s+1,m]$ and so $|2a_i-2a_h| \le 2(s-1)$, imply  
$|s_\tau(x_{ij})- s_\tau(x_{hr})| \ne 0$. 
Thus $s_\tau(u) \ne  s_\tau(v)$. 

{\bf \noindent  Case 6}: $u$ is a degree 2 vertex of $T$ and $v \in U$.  

Since $s_\tau(v) \ge p+1$, it follows that if $u\in V(P)$ or $s_\tau(u)<0$,
then $s_\tau(u) \ne s_\tau(v)$. Thus we assume $u\in V(T)\setminus V(P)$
and $s_\tau(u)>0$, where $s_\tau(u)>0$ in particular,  implies $k\ge 3$.

By~\eqref{k-even-sum} and~\eqref{k-odd-sum},  when  $s_\tau(u)>0$ and $u$ is contained in a 
big leg $L_i$ for some $i\in [1,s]$, we have 
$$
 2a_i-\left \lceil \frac{k+2}{2}\right \rceil s\ge s_\tau(u)  \ge 
 \begin{cases}
2a_i- \frac{3(k-1)}{2}s, & \text{if $k$ is odd}, \\
2a_i- \frac{3(k-2)}{2}s, & \text{if $k$ is even}. 
 \end{cases} 
$$

When $s_\tau(u)>0$ and $u$ is contained in a 
small leg $L_i$ for some $i\in [1,s]$, we have 
$$
2a_i- \frac{3(k-2)}{2}s\le s_\tau(u) \le 2a_i-\frac{k-1}{2}s. 
$$

Let $e_v \in M$ be the edge  incident with $v$, where recall $M$ is the matching defined in Step 2.1.  Then 
 by the definition of $\tau_2$, we have 
$$s_\tau(v)  >m-s+t+1+\tau_2(e_v).$$
As $\tau_2(e_v) \ge m-s+1$ and $k \ge 3$, we have $s_\tau(v)  >2m-2s+3 >2m- \lceil(k+2)/2 \rceil s \ge s_\tau(u)$ when $u$ is contained in a 
 big leg.  When $u$ is contained in a small leg $L_i$, by the definition of $\tau_2$, 
 we have $a_i<\tau_2(e_v)<m-s+t+1$.  Thus $s_\tau(v)  >m-s+t+1+a_i >2a_i-(k-1)s/2 \ge s_\tau(u)$.
  Again, this gives $s_\tau(u) \ne s_\tau(v)$. 
When $v=v_{h_t}$ is a small joint but  $s_{[\overrightarrow{P}, \tau_1]}(v)  \ge m-s$, 
we have $s_\tau(v) \ge m-s+\tau_2(e_v).$ 
 As $\tau_2(e_v) \ge m-s+1$ and $k \ge 3$, we have $s_\tau(v)  \ge 2m-2s+1 >2m-2s \ge s_\tau(u)$ when $u$ is contained in a 
 big leg.  When $u$ is contained in a small leg $L_i$, by the definition of $\tau_2$, 
we have $a_i < \tau_2(e_v)$ and $a_i-s < m-s$ (note that $a_i\ne \tau_2(e_v)$ as $L_q$ is labeled using the definition of $f$ in (I) when $s_{[\overrightarrow{P}, \tau_1]}(v)  \ge m-s$ and so is treated as a big leg). 
Thus $s_\tau(u) \le a_i+a_i-s<\tau_2(e_v)+m-s$ and so $s_\tau(u) \ne s_\tau(v)$.

When $v$ is a small joint with  $s_{[\overrightarrow{P}, \tau_1]}(v) \le p$, by the definition of $f$ in (II) and (III), we have $s_\tau(v)\le p+m-(k-2)s/2$  if $k$ is even and $s_\tau(v)\le p+m-(k-1)s$ when $k$ is odd.   
On the other hand, as $m=p+ks$, we get $s_\tau(u) \ge 2(m-s+1)-(3k-6)s/2=m+p+2-(k-2)s/2 >s_\tau(v)$ if $s_\tau(u)>0$ and $u$ is contained in a small leg or $k$ is even. 
Thus  $s_\tau(u) \ne s_\tau(v)$ if $u$ is contained in a small leg or $k$ is even.  So we assume $s_\tau(u)>0$,  $u$
is contained  in a big leg and $k$ is odd. 
Then  as $(k+3)s/2 \le ks$ (recall $k\ge 3$), we have $s_\tau(u) \ge 2(m-s+1)-(3k-3)s/2 \ge m+p+2-(k-1)s >s_\tau(v)$. 
Thus  $s_\tau(u) \ne s_\tau(v)$ . 

Lastly, we assume 
$v$ is a small joint with  $s_{[\overrightarrow{P}, \tau_1]}(v)  \ge p+1$.
By Lemma~\ref{ShanLemma}, we conclude that $v=v_{h_t}$. Also by the 
argument when $v=v_{h_t}$ is a small joint with  $s_{[\overrightarrow{P}, \tau_1]}(v)  \ge m-s$, 
we assume here that $s_{[\overrightarrow{P}, \tau_1]}(v)   \le m-s-1$. By the definition of 
$\tau_4$ in Step 2.3 and~\eqref{k-even-sum} and~\eqref{k-odd-sum}, we know that 
 $s_\tau(u)$ and  $s_\tau(v)$ have different parities and so  $s_\tau(u) \ne s_\tau(v)$. 
 
 The proof is now finished.  
 \end{proof}

\section{Open problem}

In this paper, it is shown that every subdivided caterpillar  admits an antimagic orientation. Since every bipartite antimagic graph
$G$ admits an antimagic orientation, it is natural to ask that whether subdivided caterpillars are antimagic. We propose the following conjecture. 

\begin{conj} 
Every subdivided caterpillar is antimagic. 
\end{conj}

%\bibliography{BIB}

\end{document}